\documentclass[12pt]{article}
\usepackage{amssymb}
\usepackage{amsmath}
\usepackage{mathrsfs}
\usepackage{amsfonts}
\usepackage{color}
\usepackage{vmargin}
\usepackage{amsthm}
\usepackage{mathtools}
\usepackage{a4}
\usepackage{graphicx}

\long\def\symbolfootnote[#1]#2{\begingroup%
\def\thefootnote{\fnsymbol{footnote}}\footnote[#1]{#2}\endgroup}

\setmarginsrb{40mm}{20mm}{40mm}{20mm}{10mm}{10mm}{10mm}{10mm}

\newcommand{\R}{\mathbb{R}}

\newcommand{\diam}{{\rm diam}}

\definecolor{blau}{rgb}{0.1,0.0,0.9}

\definecolor{funk}{rgb}{0.1,0.4,0.9}

\newcounter{komcounter}
\numberwithin{komcounter}{section}

\def\Xint#1{\mathchoice
   {\XXint\displaystyle\textstyle{#1}}%
   {\XXint\textstyle\scriptstyle{#1}}%
   {\XXint\scriptstyle\scriptscriptstyle{#1}}%
   {\XXint\scriptscriptstyle\scriptscriptstyle{#1}}%
   \!\int}
\def\XXint#1#2#3{{\setbox0=\hbox{$#1{#2#3}{\int}$}
     \vcenter{\hbox{$#2#3$}}\kern-.5\wd0}}

\def\dashint{\Xint-}

\theoremstyle{plain}
\newtheorem{theorem}{Theorem}[section]
\newtheorem{thm}[theorem]{Theorem}
\newtheorem{lem}{Lemma}[section]

\theoremstyle{definition}
\newtheorem{defn}{Definition}[section]

\newtheorem{rem}{\textnormal{\textbf{Remark}}}

\begin{document}

\title {Puncture repair on metric measure spaces}

%\author{
%Tomasz Adamowicz{\small{$^1$}}
%\\
%\it\small Institute of Mathematics, Polish Academy of Sciences \\
%\it\small ul. \'Sniadeckich 8, 00-656 Warsaw, Poland\/{\rm ;}
%\it\small T.Adamowicz@impan.pl
%\\
%\\
\author{
Ben Warhurst{\small{$^2$}}
\\
\it\small Institute of Mathematics,
\it\small University of Warsaw,\\
\it\small ul.Banacha 2, 02-097 Warsaw, Poland\/{\rm ;}
\it\small B.Warhurst@mimuw.edu.pl
}
%
%\date{}
\maketitle

\footnotetext[2]{The author would like to thank Mike Eastwood, Katrin F\"assler, Pekka Pankka and Tomasz Adamowicz for their helpful comments. }

\footnotetext[3]{This work was partially supported by the Simons - Foundation grant 346300 and the Polish Government MNiSW 2015-2019 matching fund.}

%\footnotetext[1]{T. Adamowicz and B. Warhurst were supported by a grant Iuventus Plus of the Ministry of Science and Higher Education of the Republic of Poland, Nr 0009/IP3/2015/73.}

%%\footnotetext[1]{T. Adamowicz was supported by a grant of National Science Center, Poland (NCN),
%%UMO-2013/09/D/ST1/03681.}

\begin{abstract} Motivated by recent interest concerning "puncture repair" in the conformal geometry of compact Riemannian manifolds (\cite{MikRod} and \cite{CFrances}), a brief exposition on generalisation to the setting of quasiconformal mappings on certain metric measure spaces is presented, as well as a brief outline on removability of porous sets.    
\end{abstract}

\section{Introduction} 
In \cite{MikRod}, the authors employ an Ahlfors--Beurling extremal length type argument, and prove the following "puncture repair" theorem: 

\smallskip

\smallskip

\noindent
{\bf Theorem}. \emph{Suppose that $M$ is a compact connected conformal manifold
 and $p \in M$. Suppose that $N$ is a compact connected conformal manifold
 and $U \subset N$ an open subset such that there exists a conformal map $f:U \to M \setminus \{p\}$.
Then $f$ extends to a conformal map of $M$ onto $N$.}

\smallskip

\smallskip
Note that by conformal manifold we mean a pair $(M,[g])$ where $[g]$ is an equivalence of Riemannian metrics where $h \sim g$ if and only if $h=\lambda^2 g$ for some function $\lambda$. If a map is conformal with respect to $g$ then it is conformal with respect to any metric equivalent to $g$, thus conformal maps between conformal manifolds are well defined. 
 
In \cite{MikRod} the authors remark that although their result is stated in terms of compact manifolds, the result is in fact local. In \cite{Zor}, Zorich shows that the theorem above holds for quasiconformal immersions between Riemannian $n$-manifolds, $n \geq 3$, when the target is simply connected.

%Theorem 2. Suppose M is a compact connected conformal manifold
%and p ∈ M. Suppose N is a compact connected conformal manifold and
%U ⊂ N an open subset such that U ∼
%= M \ {p} as conformal manifolds.
%Then this isomorphism extends to N ∼
%= M.
 
 In the modern language of geometric function theory, the extremal length has evolved into the notion of the modulus of a curve family ($=$1/extremal length), since the later behaves as an outer measure. Furthermore, the generality of the spaces to which the modulus can be used to replicate and generalise theorems of Riemannian geometry is substantial, in particular, metric measure spaces of $Q$-bounded geometry are well suited to such ends, provided one is careful regarding topology and the target is $Q$-Loewner. In this setting, generalisations of the result in \cite{MikRod} can be considered in three directions: not restricting to Riemannian metrics, considering quaisconformal mappings instead of conformal mappings, and considering underlying spaces with complicated topology, although the latter is delicate, for example on a manifold with boundary difficulties arise if the mappings do not preserve boundaries. An example of what we exclude is a ball in the closed upper half space which is centred at the origin with puncture at the origin, i.e., $\{ (x,y): y\geq 0, x^2+ y^2 < r^2 \} \setminus\{(0,0)\}$, and the map $(x,y) \to (x,y+1)$. 
   
 With the caveat that a quasiconformal mapping is assumed to preserve boundaries if they are present in the domain of the map, we have the following counterpart to Theorem 2 of \cite{MikRod}:
   
 \begin{thm}\label{puncrep} Let $(X,d,\mu)$ and $(Y,d',\nu)$ be metric measure spaces of $Q$-bounded geometry, $Q>1$, where $X$ and $Y$ are path connected. Furthermore assume that $(Y,d',\nu)$ is $Q$-Loewner. Let $x_0$ be a given point in $X$ and $U$ be an open neighborhood of $x_0$. If $U_{x_0} =U\setminus \{x_0\}$ and $f:U_{x_0} \to Y$ is a quasiconformal mapping, such that for every ball $B(x_0,r) \Subset U$, the set $C$ given by
   	$$C=\partial f( \overline{B(x_0,r)}\setminus \{x_0\} ) \setminus f( \partial B(x_0,r)) $$
   	is a continuum in $Y$,  then $f(U_{x_0})$ must be of the form $U' \setminus \{y_0\}$ where $U' \subset Y $ is an open neighbourhood of a point $y_0 \in Y $. 
\end{thm}

The result is a relatively straight forward use of the properties of the conformal modulus of a curve family and the Loewner condition, moreover the proof holds in any dimension. The $Q$-bounded geometry condition implies that the standard characterisations of quasiconformality are available and equivalent locally, while the Loewner condition on the target provides a restriction on $f$ which implies the conclusion in the theorem. Moreover, all the assumptions are fulfilled by compact Riemannian manifolds and  more (see section \ref{PoincLoew} for a discussion in the sub-Riemannian case).

Puncture repair is a specific instance of the removability problem for singularities of quasiconformal mappings, and has been a central question in geometric function theory since the pioneering work of Ahlfors and Beurling \cite{AhlfBeur}. At a general level, the removability problem seeks conditions on a compact subset $E \subset X$ such that each quasiconformal mapping of $X \setminus E$ that maps $X \setminus E$ into $Y$ extends to a quasiconformal map of $X$ onto $Y$. For example when $X$ and $Y$ are the Euclidean space $\R^n$, $n \geq 2$, it is known that $E$ is removable if the $(n - 1)$-dimensional Hausdorff measure of $E$ is zero while on the other hand, there are nonremovable Cantor sets of Hausdorff dimension $n - 1$ (see \cite{Bishop1}, \cite{Bishop2} and \cite{Vas}). These Euclidean results require analytic tools which are not available in the general setting, and so similar results are not readily at hand, although in the Riemannian case,  C. Frances \cite{CFrances} has obtained such a characterisation in the conformal Riemannian setting. A sufficient condition for removability at a general level involves the notion of a porous set. A be brief outline of the results of Balogh and Koskela \cite{BK} on porosity and removability will be given later.

%For $Q$-regular metric measure spaces, $Q > 1$, the validity of the $Q$-Poincar\'e inequality is equivalent to the Loewner condition. For example, Theorem 11.12 in \cite{hajlas} implies that a compact sub-Riemannian manifold supports a $Q$-Poincar\'e inequality where $Q$ is the Hausdorff dimension.  

%The Poincar\'e inequality holds for some non-Riemannian metrics on topological manifolds as well. In [244], Semmes provided a large class of topological manifolds equipped with metrics and measures supporting Poincar\'e inequalities.
%
%In order to state his result precisely, we recall that a metric space $(X ,d)$ is said to be linearly locally contractible if there exists a constant $C \geq 1$ such that every ball $B (x , r )$ can be contracted to a point inside the concentric ball $B(x , Cr )$.
%
%Theorem 14.2.3 (Semmes) Let $(X , d )$ be a metric space that is a topological $n$-manifold which is Ahlfors $n$-regular when equipped with the Hausdorff $n$-measure $\mathcal{H}_n$ . Assume also that $X$ is linearly locally contractible. Then $(X , d , \mathcal{H}_n )$ supports the $1$-Poincaré\,e inequality.

\section{Background on metric measure spaces}
The results discussed here are part of the elementary theory and the reader is referred to \cite{hk}, \cite{hein} and \cite{hkstB} for details.
 
Let $X=(X , d , \mu)$ be a metric measure space, that is, $(X , d )$ is separable as a metric space and $\mu$ is a locally finite Borel regular measure. Following the standard construction, one defines the length of curves defined on compact intervals and the associated arc length parameterisation which further facilitates the definition of the line integral of a real valued function on $X$, see for example \cite[chapter 5.1]{hkstB}.  If $\gamma$ is a curve defined on a noncompact interval $I$ then it's length can be defined by taking the supremum of the lengths of all curves $\gamma|_J$ where $J \subset I$ is any compact subinterval. A curve is rectifiable if it has finite length and is locally rectifiable if it is rectifiable on any compact subinterval of the domain of definition.

Let $\Gamma$ be a family of curves in $X$. The $p$-modulus of $\Gamma$  is defined as follows:
\[
{\rm Mod}_p \Gamma:=\inf_{\varrho}\int_{X}\varrho^p d\mu,
\]
where the infimum is over all nonnegative Borel functions $\varrho:X \to[0,\infty]$ such that
\[
\int_{\gamma} \varrho \,d l \geq 1,
\]
for every $\gamma\in \Gamma$ and $l$ denotes the arc length. We say that $\varrho$ is admissible for $\Gamma$. For example, in the Riemannian setting, the conformal factor of a conformal map, suitably scaled by the infimum of the set of lengths of all curves in $\Gamma$, is admissible
 
 It can also be shown that for all $p>1$, the set of nonrectifiable curves in $\Gamma$ has $p$-modulus equal to $0$, so we can assume that the curves in  $\Gamma$ are at least locally rectifiable, see \cite{hkstB} Lemma 5.3.2 and Proposition 5.3.3 page 134. 

\begin{lem}\label{lem-mod-prop}
	The following properties hold for the $p$-modulus of the curve families $\Gamma$ and $\Gamma'$:
	\begin{enumerate}
		\item If $\Gamma \subset \Gamma'$, then  ${\rm Mod}_p  \Gamma \leq {\rm Mod}_p\Gamma'$.
		\item If every curve in $\Gamma$ contains a subcurve in $\Gamma'$, then ${\rm Mod}_p  \Gamma \leq {\rm Mod}_p\Gamma'$.
	\end{enumerate}
\end{lem}
\begin{defn}\label{confMod} The "conformal modulus" is the case $p=Q$ where $Q$ is the Hausdorff dimension  of the underlying metric space. It is conformally invariant and quasi-invariant with respect to quasiconformal mapping (see Theorem \ref{eqdefs} below).     
\end{defn} 
\begin{defn}\label{genericMod}
	The notation $\Gamma(E,F,U)$ will denote the family of all curves in an open subset $U$ of $X$
	joining two disjoint closed subsets $E$ and $F$ of $U$.
\end{defn} 

With the modulus at hand we can outline the proof.  Let $B=B(x_0,r_0) \Subset U $ be an open ball of radius $r_0$ and for each $i \in \mathbb{N} $ let $B_i=B(x_0,r_i)$ be open balls where the sequence of radii $\{r_i\}$ tends to $0$.

Let $\Gamma_i$ denote the family of rectifiable curves in $\bar B \setminus B_i$ which connect $\partial B$ to $\partial B_i$. Owing to the fact that $X$ is path connected, the balls have nonempty boundary so $\Gamma_i$ is nonempty.   By Lemma \ref{upperbndmodcnd} (see below) we have ${\rm Mod}_Q\Gamma(\overline{B}_i, B \setminus B_1,B ) \to 0$ as $r_i \to 0$ and since $\Gamma_i \subset \Gamma(\overline{B}_i, B \setminus B_1,B )$ we have ${\rm Mod}_Q  \Gamma_i \to 0$ as $r_i \to 0$ by property 1 in Lemma \ref{lem-mod-prop}.

Now use the geometric properties of a quasiconformal map $f$ to conclude that the image curve family $f(\Gamma_i)$ must also have $Q$-modulus tending to zero. This last conclusion conflicts with a $Q$-Loewner condition unless $C$ is a point. 

Recall that we call a metric space $X$ endowed with a Borel measure $\mu$ is Ahlfors regular of dimension $Q$ ($Q$-regular for short) if there exists a constant $C \geq 1$ so that
\begin{align} C^{-1} r^Q \leq \mu(B, r ) \leq C r^Q \label{ahlfQreg} \end{align}
for every ball $B(x,r)  \in X$ with radius $r < {\rm diam}\,  X$. It is easy to see that if a locally
compact metric space $X$ satisfies \eqref{ahlfQreg} for some Borel measure $\mu$ then in fact it
satisfies for Hausdorff $Q$-measure (possibly with a different constant) and any measure satisfying \eqref{ahlfQreg} is comparable to the Hausdorff $Q$-measure. For this reason, the metric space $(X,d)$ is said to be $Q$-regular without specifying the measure $\mu$.

\begin{defn}\label{Qboundedgeom} (\cite{hkst}, Def 9.1). A metric measure space $(X, d, \mu)$ is said to be of locally $Q$-bounded geometry, $Q > 1$, if $(X,d)$ is separable, pathwise connected, locally compact, and if there exists constants $C_0 \geq 1$, $0<\lambda \leq 1$, and a decreasing function $\psi:(0,\infty) \to (0,\infty)$ so that the following holds: each point $x_0$ in $X$ has a relatively compact neighbourhood $U$, such that 
	\begin{align} \mu(B(x_0,r)) \leq C_0 r^Q \label{Qb1}\end{align}
	for all balls $B(x_0,r) \subset U$, and that \begin{align} {\rm Mod}_Q \Gamma(E,F,B(x_0,r)) \geq \psi(t) \label{Qb2}\end{align} whenever $B(x_0,r) \subset U$ and $E$ and $F$ are two disjoint, nondegenerate continua in $B(x_0,\lambda r)$ with ${\rm dist} (E,F) \leq t \min \{ {\rm diam} E, {\rm diam} F \}$.
\end{defn}
\begin{rem} See \cite{hkst} for more detail concerning the following consequences:	
	\begin{itemize}
		\item[(i)] Condition \eqref{Qb2}  implies \begin{align} \frac{1}{C_1} r^Q \leq \mu(B(x_0,r)), \quad B(x_0,r) \subset U,
		\end{align} for some $C_1 \geq 1$ depending only on $\psi$ and $\lambda$, and so  $(X, d, \mu)$ is uniformly locally Ahlfors $Q$-regular.
		\item[(ii)] Condition  \eqref{Qb2} can be replaced with the weaker condition \begin{align} {\rm Mod}_Q \Gamma(E,F,X) \geq \psi(t) \label{Qb3}\end{align} where $E$ and $F$ are disjoint nondegenerate continua in $B(x_0,\lambda r) \subset U$ with ${\rm dist} (E,F) \leq \min \{ {\rm diam} E, {\rm diam} F \}$. Indeed \eqref{Qb2} and  \eqref{Qb3} imply \eqref{Qb1}.  
		\item[(iii)] Connected open subsets of spaces of locally $Q$-bounded geometry are again of locally $Q$-bounded geometry, with the same data.
		\item[(iv)] Riemannian $n$-manifolds are of locally $n$-bounded geometry along with many other examples.
	\end{itemize}
\end{rem}	

A significant feature of locally $Q$-bounded geometry, stemming from item (iii) in Remark 1, is that quasiconformal mapping theory is local (in more general metric spaces, many properties of qc-maps may only hold if the map is globally defined, in particular the equivalence of conditions (i) and (ii) in the Theorem \ref{eqdefs} below.

\begin{defn}\label{loew} The statement that $(Y,d',\nu)$ is $Q$-Loewner means that there exists a decreasing function $\phi:(0,\infty) \to (0,\infty)$ so that
	\begin{align*}{\rm Mod}_Q \Gamma(E,F,Y) \geq \phi(t) \end{align*} whenever $E$ and $F$ are two disjoint, nondegenerate continua in $Y$ with ${\rm dist} (E,F) \leq t \min \{ {\rm diam} E, {\rm diam} F \}$.
\end{defn}
We note that compact Riemannian $n$-manifolds are $n$-Loewner, see \cite{hk} page 39 and \cite{hkstB} page 392.
 
\begin{thm}\label{eqdefs}(\cite{hkst}, Thm 9.8). Let $(X, d, \mu)$ and $(Y, d', \nu)$ be spaces of locally $Q$-bounded geometry, $Q > 1$ and let $U \subseteq X$ and $U' \subseteq Y$ be open neighborhoods. Then the following conditions are equivalent for a homeomorphism $f:U \to U'$:
	\begin{itemize}
		\item[(i)] there exists a constant $H \geq 1$ such that $$ \limsup_{r \to 0} \frac{ \max\{ d'(f(x),f(y)) : d(x,y) \leq r \}}{\min\{ d'(f(x),f(y)) : d(x,y) \geq r \} } \leq H;$$ 
		\item[(ii)] there exists a constant $C \geq 1$ such that
		$$ \frac{1}{C} {\rm Mod}_Q \Gamma  \leq {\rm Mod}_Q f (\Gamma ) \leq C {\rm Mod}_Q \Gamma$$
		whenever $\Gamma$ is a family of curves in $U$ and $f (\Gamma)$ denotes the family of curves
		in $U'$ obtained as images of curves from $\Gamma$ under $f$;
	\end{itemize}
\end{thm}

\begin{rem} The homeomorphism $f$ is said to be metrically quasiconformal if it satisfies (i) and geometrically quasiconformal if it satisfies (ii). The best constants $C$ and $H$ are dependent on each other and coincide when either is $1$. Thus for $1$-quasiconformal maps the quantity ${\rm Mod}_Q \Gamma$ is invariant in the sense that  $${\rm Mod}_Q f(\Gamma) ={\rm Mod}_Q \Gamma$$ for all curve families in the domain of $f$. For a direct proof in the Riemannian setting see \cite{hein} Theorem 7.10 page 51. 
\end{rem}

The extremal length method used in \cite{MikRod} considers curve families which wind around the puncture, where as the approach here will consider curve families which terminate at the puncture since the $Q$-modulus of such families is relatively well understood. In particular we have the following standard result:
 
\begin{lem}\label{upperbndmodcnd}(\cite{hkstB}, Proposition 5.3.9). Let $x_0 \in X$. If there exists a constant $C_0 > 0$ such that
$$\mu(B(x_0, r )) \leq C_0 r^Q$$
for each $0 < r < R_0$ then
$${\rm Mod}_Q \Gamma(\overline{B(x_0 , r)}, X \setminus B (x_0 , R),X) \leq C_1 \log \left ( \frac{R}{r}\right)^{1-Q}$$
whenever $0 < 2r < R < R_0$, where $C_1$ depends only on $C_0$ and $Q$.
\end{lem}
It should be note that Lemma \ref{upperbndmodcnd} does not require $Q$-bounded geometry nor a Loewner condition. Moreover, the volume growth condition is a feature of Riemannian $n$-manifolds with $Q=n$.

 \section{Proof of Theorem \ref{puncrep}} 
  
 By local $Q$-boundedness, there is a ball $B=B(x_0,r_0)  \subset U$ such that \eqref{Qb1} holds for all balls $B(x_0,r)$ such that $r<r_0$. Let $B_i=B(x_0,r_i)$ be open balls where the sequence of radii $\{r_i\}$ tends to $0$, and satisfies $2r_i <r_0$ so that Lemma \ref{upperbndmodcnd} is applicable to the curve family $\Gamma(\overline{B}_i, B \setminus B_1,B )$, that is $${\rm Mod}_Q\Gamma(\overline{B}_i, B \setminus B_1,B ) \to 0$$ as $r_i \to 0$.
  
  Let $\Gamma_i$ denote the family of rectifiable curves in $\bar B \setminus B_i$ which connect $\partial B$ to $\partial B_i$. Since $\Gamma_i \subset \Gamma(\overline{B}_i, B \setminus B_1,B )$, we have ${\rm Mod}_Q  \Gamma_i \to 0$ as $r_i \to 0$ by property 1 in Lemma \ref{lem-mod-prop}.
  
  Since $f$ is quasiconformal, Theorem \ref{eqdefs} implies that ${\rm Mod}_Q  f(\Gamma_i) \to 0$. Furthermore, by assumption, the sets $E=f( \partial B (x_0,r_0) ) $ and $F=C$, are continua in $Y$, and since every  $\gamma \in f(\Gamma_i)$ is a subcurve of $\Gamma(E, F,Y)$, property 2 of Lemma \eqref{lem-mod-prop} implies that
 $${\rm Mod}_Q\Gamma(E, F,Y) \leq  {\rm Mod}_Q f(\Gamma_i).$$ We note that $\diam E >0$ and so $\diam F$ must be zero, for if it were otherwise, the Loewner condition on $Y$ would imply that ${\rm Mod}_Q\Gamma(E, F,Y)$ is nonzero thus prohibiting ${\rm Mod}_Q  f(\Gamma_i) \to 0$.
 
 \section{Poincar\'e inequalities}\label{PoincLoew}
Validating the Loewner condition for a metric measure space is a matter typically addressed via Poincar\'e inequalities. More precisely, if $(X,d,\mu)$ is proper (closed balls are compact), satisfies a convexity condition (see \cite{hk} section 2.15), and $Q$-regular, $Q > 1$, then the $Q$-Loewner condition is equivalent to the validity of the $Q$-Poincar\'e inequality, see \cite{hk} Corollary 5.13. 
%Under suitable (weak) topological conditions, a Q -regular space X is Loewner if and only if it supports a (1, Q) -Poincaré inequality [HK3, Corol-
%lary 5.12].

In the setting of metric measure spaces with no Riemannian structure, the following notion of upper gradients, formulated by Heinonen and Koskela in \cite{hk},
plays the role of derivatives. A Borel function $g$ on $X$ is an upper gradient of a
real-valued function $u$ on $X$ if for all nonconstant rectifiable paths $\gamma: [0, l_\gamma ] \to X$
parameterized by arc length,
$$ |u(\gamma(0)) - u(\gamma(l_\gamma))| \leq \int_{\gamma } g \, d l$$
where the above inequality is interpreted as saying also that $\int_{\gamma } g \, d l=\infty$ whenever
$|u(\gamma(0))|$ or  $|u(\gamma(l_\gamma))|$ is infinite. %If the above inequality fails only for a
%curve family with zero $p$-modulus, then $g$ is a $p$-weak upper gradient of $u$. It is
%known that the $L^p$-closed convex hull of the set of all upper gradients of $u$ that
%are in $L^p(X)$ is precisely the set of all $p$-weak upper gradients of $u$ in $L^p(X)$, see
%Lemma 2.4 in \cite{KoMc}.

\begin{defn}\label{poincareineq} 
	We say that $(X, d, \mu)$  supports a $p$-Poincar\'e inequality if there are constants $\tau, C > 0$ such that for all balls $B \subset X$, all measurable functions $f$
on $X$, and all $p$-weak upper gradients $g$ of $f$,
$$ \dashint_B |f - f_B| d \mu \leq C r \left (
\dashint_{τB} g^p d \mu \right )^{1/p}$$
where $r$ is the radius of $B$ and
$$f_B =\dashint_B f\,  d \mu =\frac{1}{\mu(B)} \int_B f \, d\mu. $$
\end{defn}

%Recall that we call a metric space X endowed with a Borel measure μ an
%Ahlfors regular space of dimension Q (for short, a Q -regular space) if there exists
%a constant C 0 ≥ 1 so that
%(1.3)
%C 0 −1 r Q ≤ μ(B r ) ≤ C 0 r Q
%for every ball B r in X with radius r < diam X . It is easy to see that if a locally
%compact metric space X satisfies (1.3) for some Borel measure μ then in fact it
%satisfies it for Hausdorff Q -measure H Q (possibly with a different constant C 0 )
%and any measure satisfying (1.3) is comparable to H Q . For this reason, we will
%often say that a metric space X is Q -regular without specifying the measure in
%question. Q -regularity of a metric space X can be regarded as a quantitative
%and scale-invariant version of the qualitative statement that X has Hausdorff
%dimension Q .

\begin{thm}\label{Qboundedgeom} (\cite{Shan}, Def 2.2) A metric measure space $(X, d, \mu)$ is of locally $Q$-bounded geometry, $Q > 1$, if and only if $\mu$ is (locally) Ahlfors $Q$-regular and the space supports a local $p$-Poincar\'e inequality for some $1 \leq p < Q$ in the sense of Definition \ref{poincareineq}.
\end{thm}

 For example, Carnot groups with their sub-Riemannian metric and Lebesgue measure, support a $1$-Poincar\'e inequality and, by H\"older's inequality, a $Q$-Poincar\'e inequality also holds. Thus Carnot groups have locally $Q$-bounded geometry and are $Q$-Loewner. In the setting of a set of equiregular H\"ormander vector fields on an open set $\Omega \subseteq \R^n$, together with the corresponding Carnot-Caratheodory distance and  Lebesgue measure, Jerison \cite{Jerison} shows that a $1$-Poincar\'e inequality holds. It follows that compact equiregular Carnot-Caratheodory/sub-Riemannian manifolds support a $1$-Poincar\'e inequality and thus have locally $Q$-bounded geometry and are $Q$-Loewner. Hence it follows that Theorem \ref{puncrep} applies in this setting. 

\section{Porous sets}

At the metric measure space level with "separation axiom {\bf S}", Balogh and Koskela \cite{BK} consider conditions on a compact subset $E$ of a $Q$-regular Loewner space that guarantee that each quasiconformal mapping of $X \setminus E$ that maps $X \setminus E$ into $X$ and maps bounded sets to bounded sets, extends to a quasisymmetric homeomorphism of $X$ onto $X$ (note that when a homeomorphism maps bounded sets to bounded sets then quasiconformality and quasisymmetry are equivalent conditions). The condition they arrive at (see Theorem 1.1 \cite{BK}) is called spherical $t$-porosity of $E$, that is to say: there exits $t>1$ such that for each $x \in E$ there is a sequence $r_j$ tending to zero with $E \cap B(x, t r_j ) \setminus  B(x, t^{-1}  r_j ) = \emptyset$ for all $j$. 

Strictly speaking Theorem 1.1 in \cite{BK} assumes $X$ is unbounded and that $X$ is also the target, however in the remarks on page 573 of \cite{BK}, it is explained that these assumptions can be removed. In particular, for a general target $Y$, their result holds when $X$ is proper, $Q$-regular with $Q > 1$, and Loewner, while the target space $Y$ is assumed to be proper, $Q$-regular with the same $Q$, $C$-linearly locally connected, and satisfies the separation axiom {\bf S}. In particular, Theorem A.2 in \cite{BK} shows that metric connected $n$-manifolds, $n \geq 2$, satisfy separation axiom {\bf S}. Consequently the puncture repair theorem holds for porous $E$ in the Riemannian setting.

{\bf Axiom S}. A metric space $X$ satisfies the separation axiom S if, for each compact subset $E$ of $X$, there is a constant $\lambda > 0$ so that the following holds: If a finite union $A \subset F$ of pairwise disjoint continua with diameters less than $\lambda$ separates points $x_1, x_2 \in X$, then a component of $A$ separates $x_1$ and $x_2$.


\begin{thebibliography}{99}

\bibitem{AhlfBeur} {\sc L. Ahlfors and A. Beurling}, \emph{Conformal invariants and function-theoretic null-sets}, Acta Math. 83 (1950), 101--129.

\bibitem{BK} {\sc Z. Balogh and P. Koskela}, \emph{Quasiconformality, quasisymmetry, and removability in Loewner spaces}, Duke Math. J., 101(3):554--577, (2000).

\bibitem{Bishop1} {\sc C. Bishop}, \emph{Some homeomorphisms of the sphere conformal off a curve}, Ann. Acad. Sci. Fenn. Ser. A I Math. 19 (1994), 323--338.

\bibitem{Bishop2}  {\sc C. Bishop}, \emph{Non-removable sets for quasiconformal and bilipschitz mappings in $\R^n$}, arXiv:math/9806015, (1998).

\bibitem{MikRod} {\sc M. Eastwood, R. Gover}, \emph{Volume growth and puncture repair in conformal geometry}, arXiv:1707.03087, (2017).

\bibitem{CFrances} {\sc C. Frances}, \emph{Removable and essential singular sets for higher dimensional conformal maps}, Comment. Math. Helv. 89 (2014) 405--441.

\bibitem{hajlas} { \sc P. Hajłasz and P. Koskela},  \emph{Sobolev met Poincar\'e}, Mem. Amer. Math. Soc., (2000), 145(688).

\bibitem{hein} {\sc J. Heinonen}, \emph{Lectures on Analysis on Metric Spaces}, Springer, New York, NY, (2001).

\bibitem{hk} {\sc J. Heinonen, P. Koskela}, \emph{Quasiconformal maps in metric spaces with controlled geometry}, Acta Math. 181 (1998), 1--61.

\bibitem{hkst} {\sc J. Heinonen, P. Koskela, N. Shanmugalingam and J. Tyson}, \emph{Sobolev classes of Banach space-valued functions and quasiconformal mappings}, J. Anal. Math. 85 (2001), 87--139.

\bibitem{hkstB} {\sc J. Heinonen, P. Koskela, N. Shanmugalingam and J. Tyson}, \emph{Sobolev Spaces on Metric Measure Spaces: An Approach Based on Upper Gradients}, New Mathematical Monographs:27, Cambridge University Press, Cambridge United Kingdom, (2015).

\bibitem{Jerison} {\sc D. Jerison}, \emph{The Poincar\'e inequality for vector fields satisfying H\"ormander’s condition}, Duke Math. J., Volume 53, Number 2 (1986), 503-523.


\bibitem{Shan} {\sc N. Shanmugalingam}, \emph{Singular behavior of conformal Martin kernels, and nonitangential limits of conformal mappings}, Ann. Acad. Sci. Fenn. Ser. A I Math.,
Volume 29, 2004, 195--210.

\bibitem{Vas} {\sc J. V\"ais\"al\"a}, \emph{Removable sets for quasiconformal mappings}, J. Math. Mech. 19 (1969/1970), 49--51.

\bibitem{Zor} {\sc V. A. Zorich}, \emph{A nonremovable singularity of a quasiconformal immersion}. Russ. Math. Surv. 64 (2009), 173--174.
\end{thebibliography}
\end{document}